\documentclass[12pt,a4paper]{amsart}

\usepackage[english]{babel}
\usepackage{vmargin}
\usepackage[dvips]{graphicx}
\usepackage{color}
\usepackage{url}
\usepackage{microtype}

\usepackage{amssymb}
\input{epsf}
\usepackage{epsfig}
\usepackage{a4wide}
\usepackage{supertabular}
\usepackage{array}
\usepackage{multicol}
\usepackage[utf8x]{inputenc}
\usepackage[T1]{fontenc}

\setmarginsrb{3cm}{2cm}{3cm}{3cm}{75pt}{20pt}{20pt}{30mm}
\renewcommand{\phi}{\varphi}
\renewcommand{\geq}{\geqslant}
\renewcommand{\leq}{\leqslant}

\def\R{\mathbb R}
\def\N{\mathbb N}

\def\E{\mathbb E}

\def\eps{\varepsilon}

\def\eps{\varepsilon}

\def\tilde{\widetilde}

\def\vol{{\rm vol}\,}

\def\<{\langle}
\def\>{\rangle}

\newcommand{\abs}[1]{\left\vert#1\right\vert}

\def\={\mathrel{\mathop:}=}
\newcommand{\eqdef}{=\mathrel{\mathop:}}

\renewcommand{\d}{\mathrm{d}\hspace{-0.02em}}

\def\begeq{\begin{equation}}
\def\endeq{\end{equation}}

\newtheorem{Thm}{Theorem}
\newtheorem{Lem}[Thm]{Lemma}

\newtheorem{Cor}[Thm]{Corollary}
\newtheorem{Prop}[Thm]{Proposition}

\numberwithin{equation}{section}
\theoremstyle{definition}

\theoremstyle{remark}

\newtheorem*{Thm*}{Theorem}
\newtheorem*{Lem*}{Lemma}
\newtheorem*{Conj*}{Conjecture}
\newtheorem*{Cor*}{Corollary}
\newtheorem*{Def*}{Definition}
\newtheorem*{Prop*}{Proposition}
\newtheorem*{Exo*}{Exercise}
\newtheorem*{Exs*}{Examples}
\newtheorem*{Ex*}{Example}
\newtheorem*{Rk*}{Remark}
\newtheorem*{Rks*}{Remarks}

\DeclareMathOperator*{\Supp}{Supp}

\def\Ack{\medskip\noindent {\bf Acknowledgements:}\ \ignorespaces}

\begin{document}

\title[A curved Brunn--Minkowski inequality on the discrete hypercube]
{A curved Brunn--Minkowski inequality on the discrete hypercube
\\
\footnotesize
Or: What
is the Ricci curvature of the discrete hypercube?}
%\date{\today}

\author{Y. Ollivier}
\author{C. Villani}

\begin{abstract}
We compare two approaches to Ricci curvature on non-smooth spaces, in the
case of the discrete hypercube $\{0,1\}^N$. While the coarse Ricci
curvature of the first author readily yields a positive value for curvature, the
displacement convexity property of Lott, Sturm and the second author could not be
fully implemented. Yet along the way we get new results of a
combinatorial and probabilistic nature, including a curved
Brunn--Minkowski inequality on the discrete hypercube.
\end{abstract}

\maketitle

\setcounter{section}{0}
%\tableofcontents

\section*{Introduction}

%TODO:
%Entropy: case where mu_0 close to mu_1: HWI inequality?

Let $A_0, A_1$ be two compact, nonempty subsets of $\R^n$. In one of its guises,
the remarkable Brunn--Minkowski inequality states that
\[
\ln \vol A_t \geq (1-t) \ln \vol A_0+t\ln \vol A_1
\]
where $0\leq t\leq 1$ and
$A_t=\{(1-t)a_0+ta_1, \;a_0\in A_0, \, a_1\in A_1\}$ is the set of
$t$-midpoints between $A_0$ and $A_1$. In other words, the logarithm of
the volume of $A_t$ is concave. We refer to \cite{Gar02} for a nice
survey. This is the ``infinite-dimensional''
version of the Brunn--Minkowski inequality, from which the more common
version using $1/n$-th powers instead of logarithms can be derived (see
Eq.~(22) in \cite{Gar02}).

If $\R^n$ is replaced with a Riemannian manifold, the presence of
positive curvature \emph{improves} this inequality. Indeed, in
\cite{CMS06} (elaborating on \cite{CMS01}) it is
proved that if $X$ is a smooth and complete Riemannian manifold with Ricci
curvature at least $K$ for some $K\in \R$, then for any two compact,
nonempty
%%CMS says nothing about compact, but it seems to be necessary in the
%Euclidean case...
subsets $A_0, A_1\subset X$, we have
\[
\ln \vol A_t \geq (1-t) \ln \vol A_0+t\ln \vol
A_1+\frac{K}{2}\,t(1-t)\,d(A_0,A_1)^2.
\]
Here the set of $t$-midpoints $A_t$ is defined as the set of all
$\gamma(t)$ where $\gamma$ is any minimizing geodesic such that
$\gamma(0)\in A_0$ and $\gamma(1)\in A_1$. The distance $d(A_0,A_1)$ is
$\inf_{a_0\in A_0, \,a_1\in A_1} d(a_0,a_1)$.

\begin{figure}
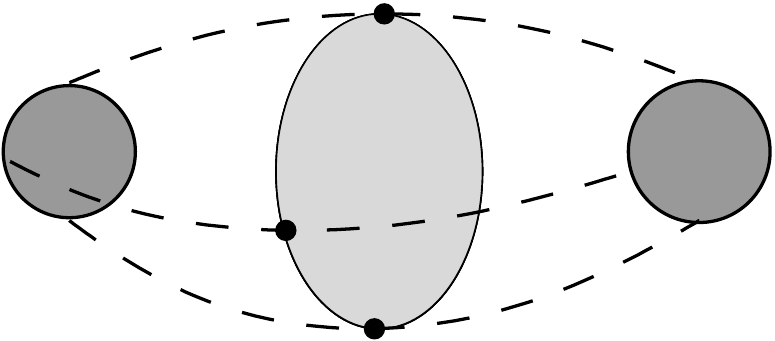
\caption{In positive curvature, midpoints spread out.}
\label{fig1}
\end{figure}

Actually this kind of inequality has been used as a tentative
\emph{definition} of positive Ricci curvature on more general, non-smooth
spaces. The idea is that, in positive curvature, ``midpoints spread out''
so that the set of midpoints of two given sets is larger than in the
reference Euclidean case (Fig.~\ref{fig1}). This led to the notion of
\emph{displacement convexity of entropy} for Riemannian manifolds
\cite{RS05,CMS01,OV00}, later developed by Sturm \cite{Stu06} and Lott and the
second author \cite{LV09}. However, it is not clear how this fares for
discrete spaces \cite{BS09}.

\begin{figure}
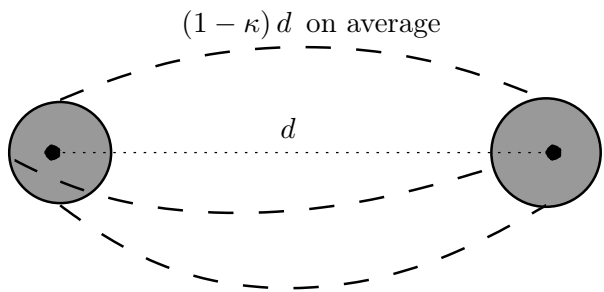
\caption{In positive curvature, balls are closer than their centers.}
\label{fig2}
\end{figure}

Another approach to define the Ricci curvature of discrete spaces is
\emph{coarse Ricci curvature}, developed by the first author
\cite{curvmarkov_cras,curvmarkov}. The motto is that, in
positive curvature, ``balls are closer than their centers are'' in
transportation distance (Fig.~\ref{fig2}).

We compare both approaches applied to the discrete hypercube
$X=\{0,1\}^N$. This is the most simple discrete space expected to have
positive Ricci curvature in some sense, for a variety of reasons (see,
e.g., paragraph $3\tfrac12$.21 ``Spheres, cubes, and the law of large
numbers'' in \cite{Gro99}).
The subtitle question ``What is the Ricci curvature of the discrete
hypercube?'' was asked verbatim by Stroock in a seminar as early as 1998,
in a context of logarithmic Sobolev inequalities.

The formalism of coarse Ricci curvature is readily available for the
hypercube
and yields a
value of $\frac{2}{N+1}$ for the Ricci curvature of $\{0,1\}^N$
(section~\ref{sec:coarsericci}).
On the other hand, we could not fully implement the displacement
convexity of entropy (properly discretized) in the hypercube. Yet, along
the way, we still get a combinatorial \emph{Brunn--Minkowski inequality on the
hypercube}, including a positive curvature term. The resulting value of
curvature is $\approx 1/N$, compatible with coarse Ricci curvature.

\Ack The authors would like to thank Prasad Tetali for helpful comments
on concentration in the symmetric group, which led to improved constants.

\section{Statement of results}

\subsection{Brunn--Minkowski inequality in the hypercube.}

We consider the discrete hypercube $X\=\{0,1\}^N$, $N\in\N$,
equipped with the Hamming (or $\ell^1$) metric
\[
d((x_i),(y_i))\=\#\{i,\;x_i\neq y_i\}.
\]
For $A$ and $B$ nonempty subsets of $X$, we define
$
d(A,B)\=\inf_{a\in A,b\in B} d(a,b).
$

Let $a$ and $b$ be two points in $X$. A
\emph{midpoint} of $a$ and $b$ is any point $m$ such that
$d(m,a)+d(m,b)=d(a,b)$ and $\abs{d(m,a)-d(a,b)/2}<1$. More explicitly: if
$d(a,b)$ is even, a midpoint is the middle point on any shortest path
from $a$ to $b$ in $X$, and if $d(a,b)$ is odd, a midpoint is one the two middlemost
points on such a shortest path.
In the hypercube, midpoints are by no means unique: the number of midpoints of $a$ and $b$
is the binomial coefficient $\binom{d(a,b)}{d(a,b)/2}$ if $d(a,b)$ is even, and
$2\binom{d(a,b)}{(d(a,b)-1)/2}$ if $d(a,b)$ is odd.

If $A$ and $B$ are two subsets of $X$, the set of midpoints of $A$ and
$B$ is the set of midpoints of all pairs $(a,b)\in A\times B$.

\begin{Thm}
\label{thm:BMcube}
Let $A$ and $B$ be two nonempty subsets of $\{0,1\}^N$. Let $M$ be the
set of midpoints of $A$ and $B$. Then
\[
\ln \#M \geq \frac12 \ln \#A+\frac12\ln \#B + \frac{K}{8}\,d(A,B)^2
\]
with $K=\frac{1}{2N}$.
\end{Thm}

This is analogous to the curved Brunn--Minkowski inequality
above in Riemannian manifolds (for $t=1/2$), with $K$ playing the role of
a curvature lower bound.

The order of magnitude $\frac1N$ for $K$ is optimal: indeed, when $A$ and
$B$ are singletons lying at distance $N$, then $d(A,B)^2=N^2$, while the
number of midpoints is $\binom{N}{N/2}\sim 2^N\sqrt{\frac{2}{\pi N}}$, so that
$\ln \#M$ grows linearly in $N$.

We will now see that this theorem can be improved by replacing $d(A,B)$ with a
transportation distance.

\subsection{Entropy of midpoints in the hypercube.}
Theorem~\ref{thm:BMcube} appears as a particular case of a refined
statement using probability measures instead of sets.

Let $\mu$ be a probability measure on a discrete set $X$. Its
\emph{Shannon entropy} is
\[
S(\mu):=-\sum_{x\in X} \mu(x) \ln \mu(x).
\]
In particular, if $\mu$ is the uniform distribution on a finite subset
$A\subset X$, then $S(\mu)=\ln\#A$.

In this paper, we shall also use the \emph{relative entropy} (or
Kullback--Leibler divergence) of a measure $\mu$ with respect to a
reference probability measure $\nu$, defined as
\[
H(\mu|\nu):=\sum_{x\in X} \mu(x) \ln\frac{\mu(x)}{\nu(x)}\geq 0.
\]
If $X$ is finite and the reference measure $\nu$ is uniform on $X$, then
we have $H(\mu|\nu)=\ln\#X-S(\mu)$.

To state an entropic version of Theorem~\ref{thm:BMcube} we define the
midpoints of two measures as follows. Loosely speaking,
we first pick a random point $a$ under $\mu_0$, then an
\emph{independent} random point $b$ under $\mu_1$, and finally we pick a
random midpoint of $a$ and $b$ uniformly over all such midpoints.

\newcommand{\midm}{\mathit{mid}}

More precisely, let $a$ and $b$ be two points of the hypercube $X$. The \emph{midpoint
measure} $\midm(a,b)$ is defined as the uniform probability measure on all
midpoints of $a$ and $b$.
Let now $\mu_0, \mu_1$ be two probability measures on $X$. The
\emph{midpoint measure} of $\mu_0$ and $\mu_1$ is defined as
\[
\midm(\mu_0,\mu_1)\= \iint \midm(a,b)\,\d\mu_0(a)\d\mu_1(b).
\]

\begin{Thm}
\label{thm:entropicBM}
Let $\mu_0$ and $\mu_1$ be two probability measures on the discrete hypercube
$X=\{0,1\}^N$. Let $\mu_{1/2}=\midm(\mu_0,\mu_1)$ be their midpoint
measure. Then
\[
S(\mu_{1/2})\geq\frac12
\left(S(\mu_0)+S(\mu_1)\right)+\frac{K}{8}\,W_1(\mu_0,\mu_1)^2
\]
with $K=\frac1{2N}$.
Equivalently,
\[
H(\mu_{1/2}|\nu)\leq \frac12
\left(S(\mu_0|\nu)+S(\mu_1|\nu)\right)-\frac{K}{8}\,W_1(\mu_0,\mu_1)^2
\]
with $\nu$ the uniform probability measure on $\{0,1\}^N$.
\end{Thm}

Here we use the
\emph{$L^1$ Wasserstein distance}
\[
W_1(\mu,\mu')\=
\inf_{\xi}\iint d(a,b)\,\d\xi(a,b)
\]
where the infimum is taken over all measures $\xi$ on $X\times X$ such
that $\int_b \d\xi(a,b)=\d\mu(a)$ and $\int_a \d\xi(a,b)=\d\mu'(b)$,
i.e., all couplings of $\mu$ and $\mu'$. We refer to \cite{Vil03} for
more background on this topic.

Note that $W_1(\mu_0,\mu_1)$ is always at least $d(A,B)$ for
$\mu_0$ and $\mu_1$ supported in sets $A$ and $B$; in particular, if
$\mu_0$ and $\mu_1$ are taken uniform in $A$ and $B$,
Theorem~\ref{thm:entropicBM} is really a refinement
of Theorem~\ref{thm:BMcube}.

\subsection{Limitations and open questions.} A first limitation of these
results is the necessity to take $t=1/2$. This comes from the combinatorial nature of our
proof, which, for the most basic situation $K=0$, consists in building an
injection from $A\times B$ into $M\times M$.

This can probably be circumvented if we assume that the sets $A$ and $B$
are convex (i.e.\ the midpoint of two points in $A$ lies in $A$, and
likewise for $B$): then, we can describe $t$-midpoints of $A$ and $B$ as
iterated $1/2$-midpoints. (If $A$ or $B$ are not convex, iterating only
yields midpoints of several points in $A$ and several points in $B$,
which is not what we want.)

The injection  from $A\times B$ into $M\times M$ used in our proof very
naturally extends to an injection from $A\times B$ into $M_t\times
M_{(1-t)}$, with $M_t$ the set of $t$-midpoints. This leads to a lower
bound for $\ln \# M_t+\ln\# M_{(1-t)}$ in terms of $\ln\#A +\ln\#B$ plus
a curvature term. This also holds in the Riemannian case (by adding
the Brunn--Minkowski inequality for $t$ and for
$(1-t)$). We do not know if there is a particular interpretation of this
inequality.

% Another very natural extension of our construction is to increase the
% number of sets involved: if three sets $A,B,C$ are given, the same proof
% yields an injection from $A\times B\times C$ to $M\times M\times M$ with
% $M$ the set of ``$1/3$-barycenters'' of $A$, $B$ and $C$. In $\R^n$,
% this Brunn--Minkowski inequality with more sets follows from the usual
% case; we do not know if something similar holds in the Riemannian case,
% or if the combinatorial version in the cube has any particular interest.

Our initial goal was to prove that the discrete hypercube has positive Ricci
curvature in the sense of Lott, Sturm and the second author, i.e., that
the hypercube satisfies displacement convexity of entropy (see below).  The main
difference with %the Brunn--Minkowski inequality is that in our result we
our result is that, in the Brunn--Minkowski inequality, we
consider all midpoints of all pairs of points $(a,b)$ with law
$\mu_0\otimes \mu_1$; whereas for displacement convexity, one should
first choose an optimal coupling between $\mu_0$ and $\mu_1$ and then
only consider the midpoints of those pairs $(a,b)$ that make up the
optimal coupling.
%For instance, if $\mu_0=\mu_1$ these are completely
%different problems.
The two properties coincide only when $\mu_0$ is
a Dirac measure, in which case our result is related to Sturm's
\emph{measure contraction property} \cite{Stu06}.

So as far as we know, the problem of computing the Ricci curvature of the
hypercube using the displacement convexity approach is still open.

\section{Two approaches to discrete Ricci curvature}

We now present in more detail the two known approaches for Ricci
curvature on discrete spaces. This is not necessary to understand our
results and proofs, but provides the original motivation.

\subsection{Coarse Ricci curvature (after the first author).}

\label{sec:coarsericci}

The basic idea of coarse Ricci curvature is to take two small balls and
compute the transportation distance between them. If this distance is
smaller than the distance between the centers of the balls, then coarse
Ricci curvature is positive.

This is formalized as follows \cite{curvmarkov_cras,curvmarkov}. Let $(X,d)$ be a metric space equipped with a
measure $\mu$. Let $\eps$ be a discretization parameter (we take $\eps=1$ for a
graph) and assume that all $\eps$-balls in $X$ have finite and non-zero
measure. For $x\in X$ define the measure $\mu_x$ by restricting $\mu$ to
the closed $\eps$-ball around $x$:
\[
\mu_x\=\frac{\mu_{|B(x,\eps)}}{\mu(B(x,\eps))}
\]
with $B(x,\eps)=\{y\in X,d(x,y)\leq \eps\}$.

If $x$ and $y$ are two points in $X$, then the \emph{coarse Ricci
curvature along $(x,y)$} is the number $\kappa(x,y)$ defined by
\[
W_1(\mu_x,\mu_y) \eqdef (1-\kappa(x,y))\,d(x,y)
\]
where $W_1$ is the $L^1$ Wasserstein distance as defined earlier. If this
is applied to a Riemannian manifold, this gives back the ordinary Ricci
curvature when $\eps\to 0$, up to scaling by $\eps^2$.

\begin{figure}
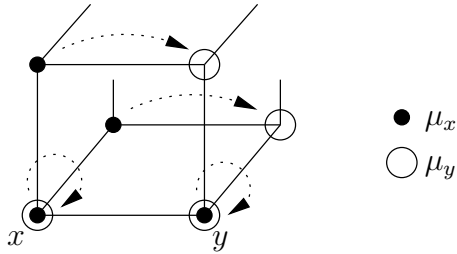
\caption{Coarse Ricci curvature in the hypercube.}
\label{fig:riccicube}
\end{figure}

Let us apply this to the discrete hypercube $X=\{0,1\}^N$ equipped with
the uniform measure. The measure $\mu_x$ is uniform on the $N+1$
neighbors of $x$ (counting $x$ itself). When $x$ and $y$ are neighbors,
it is very easy to compute the curvature $\kappa(x,y)$, as illustrated on
Figure~\ref{fig:riccicube}. Indeed, we have to move the $N+1$ neighbors
of $x$ to the $N+1$ neighbors of $y$; out of these $N+1$ points, two are
already in place ($x$ and $y$ themselves) and do not need to move, and
the others have to move by a distance $1$. So
$W_1(\mu_x,\mu_y)=1-2/(N+1)$ and $\kappa(x,y)=2/(N+1)$.

If $x$ and $y$ are not neighbors, we use a locality property of coarse
Ricci curvature. Namely, if the space $X$ is $\delta$-geodesic (i.e.\ if
the distance between two points is realized by a sequence of points with
jumps at most $\delta$), then it is enough to compute $\kappa(x,y)$ for
$d(x,y)\leq \delta$ (Exercise~2 in \cite{curvmarkov_cras}). A graph is $1$-geodesic by definition of the graph
metric, so it is enough to work with neighbors.

A lower bound on coarse Ricci curvature comes with a number of
consequences \cite{curvmarkov}. For the discrete hypercube equipped with
the uniform measure these properties were already known (but not on the
hypercube with e.g.\ $\text{Bernoulli}(\theta/N)$ measures \cite{JO10}).

In general, one may directly choose an arbitrary Markov kernel
$\mu_x$ (without using a global measure $\mu$); this leads to
interesting applications \cite{JO10}.

\subsection{Displacement convexity (after Lott, Sturm and the second
author).} In
\cite{RS05} (following ideas from \cite{OV00}), Renesse and Sturm
present a characterization of Ricci curvature on
Riemannian manifolds, based on the idea that in positive curvature,
``midpoints spread out''.

Let $X$ be a smooth, complete Riemannian manifold.
Let $\d x$ be the Riemannian volume measure on $X$. Given a probability
measure $\mu$ on $X$, define its relative entropy as $H(\mu|\d
x)\=\int \ln \frac{\d\mu}{\d x}\,\d \mu$ if the
integral makes sense, or $+\infty$ otherwise.

Let $\mathcal{P}^2(X)$ be the set of probability measures on $X$ with
finite second moment, i.e.\ those probability measures $\mu$ such that
$\int d(\mathrm{pt},x)^2\,\d\mu(x)<\infty$ for some (hence any) point
$\mathrm{pt}\in X$. On $\mathcal{P}^2(X)$, the Wasserstein distance $W_2$
is well-defined. Moreover, $\mathcal{P}^2(X)$ equipped with the metric
$W_2$ is a geodesic space: given any two probability measures
$\mu_0,\mu_1\in \mathcal{P}^2(X)$, there exists a curve
$(\mu_t)_{t\in(0;1)}$ in $\mathcal{P}^2(X)$ with
$W_2(\mu_t,\mu_{t'})=\abs{t-t'}\,W_2(\mu_0,\mu_1)$ for $t,t'\in[0;1]$.
Such a curve is called a \emph{displacement interpolation} between
$\mu_0$ and $\mu_1$. We refer to Chapter~7 of \cite{Vil08} for more
details.

Theorem~1.1 in \cite{RS05} asserts that the Riemannian manifold $X$ has
Ricci curvature at least $K\in\R$ if and only if the following
inequality is satisfied: for any two measures
$\mu_0,\mu_1\in\mathcal{P}^2(X)$, for any $W_2$-geodesic
$(\mu_t)_{t\in(0;1)}$ joining them, we have
\[
H(\mu_t|\d x)\leq (1-t)H(\mu_0|\d x)+tH(\mu_1|\d
x)-\frac{K}{2}t(1-t)W_2(\mu_0,\mu_1)^2,
\]
a property called \emph{displacement convexity} of the entropy
function.

For any probability measure $\mu$ we have $H(\mu|\d x)\geq -\ln\vol
\Supp(\mu)$, with equality when $\mu$ is uniform on its support.  Taking
$\mu_0$ and $\mu_1$ to be uniform probability distributions on sets $A_0$
and $A_1$ respectively, we see that displacement convexity of entropy
implies an inequality between the logarithms of the volumes of the
support of $\mu_t$, $\mu_0$ and $\mu_1$.  This inequality is very similar
to the Brunn--Minkowski inequality mentioned earlier. Actually, an
important property of displacement interpolation is that the measure
$\mu_t$ will charge only $t$-midpoints between the supports of $\mu_0$
and $\mu_1$ (Corollary~7.22 in \cite{Vil08}, basically due to Brenier and
McCann),  and so the Brunn--Minkowski inequality in a Riemannian manifold
really follows from convexity of entropy.

Displacement convexity of entropy makes sense in an arbitrary geodesic
space. In \cite{Stu06,LV09}, it is taken as the basis for a
notion of Ricci curvature in such spaces. The definition depends on two
parameters $K$ (the curvature) and $N$ (a ``dimension''). Displacement
convexity of entropy as written here corresponds to $N=\infty$,
the simplest and weakest case.

Interestingly, this approach applies to spaces with positive curvature in
the sense of Alexandrov \cite{Pet}.

Application to discrete spaces requires some changes: for instance, in
the case of the hypercube considered in this article, clearly if two
points are at odd distance they do not have an exact midpoint, but they
have an approximate midpoint up to an error term $\pm 1/2$.  Such an
approach is used in \cite{Bon09} to define the Brunn--Minkowski
inequality on discrete spaces. In \cite{BS09}, Bonciocat and Sturm use
approximate midpoints in the space of probability measures to extend the
definition of displacement convexity of entropy to discrete spaces, and
provide examples of planar graphs satisfying this property. To our
knowledge, these planar graphs are the only discrete examples so far.

%[Discuss relationship with our result and further work to prove it for the
%hypercube...]

\section{Brunn--Minkowski inequality without curvature}

To make the idea clearer and introduce necessary concepts, we begin with
a simplified version of Theorem~\ref{thm:BMcube}, namely
the same statement with $K=0$.
So let $A,B$ be two nonempty subsets of the hypercube $X=\{0,1\}^N$.
Let $M$ be the set of midpoints of $A$ and $B$. We want to prove that
\[
\ln \#M\geq
\frac12\left(\ln\#{A}+\ln\#{B}\right)
\]
or equivalently
\[
\#M\geq \sqrt{\#A\;\#B}.
\]

Let $a=(a_i)_{1\leq i\leq N}\in A$ and $b=(b_i)_{1\leq i\leq N}\in B$. A midpoint $m=(m_i)$ of $a$ and $b$ is a sequence of
bits such that $m_i=a_i$ whenever $a_i=b_i$ and such that half the
remaining bits coincide with those of $a$ and the other half with those
of $b$. Let $r=d(a,b)$ be the number of distinct bits between $a$ and
$b$. For fixed $a$ and $b$, there is a one-to-one correspondence between
the midpoints $m$ of $a$ and $b$
and the subsets $c\subset\{1,\ldots,r\}$
with cardinality $r/2$ (if $r$ is even) or $r/2\pm1/2$ ($r$ odd): 
among the $r$ distinct bits between $a$ and $b$, the set $c$ describes
those picked from $a$ in the construction of $m$.

We shall call \emph{$r$-crossover} such a $c\subset\{1,\ldots,r\}$ with
$\abs{\#c-r/2}\leq 1/2$. We shall denote $m=\phi_c(a,b)$ the midpoint of
$a$ and $b$ defined by crossover $c$. If $c$ is a crossover, we shall
denote by $\bar{c}$ its complement, which is also a crossover.

Note that, given a fixed $d(a,b)$-crossover $c$, the pair
$\Phi_c(a,b)\=(\phi_c(a,b),\phi_{\bar{c}}(a,b))=(m,m')$ allows to recover
$a$ and $b$.  Indeed, the identical bits in $m$ and $m'$ are the same as
in $a$ and $b$; the bits that differ between $m$ and $m'$ also differ
between $a$ and $b$, and knowledge of the crossover $c$ tells us exactly
which of those come from $a$ or $b$.

In particular, for each $r\in\{0,\ldots,N\}$, let us define the
$r$-crossover $c_r\=\{1,2,\ldots,\lfloor r/2\rfloor\}$.
Then the map $(a,b)\to \Phi_{c_{d(a,b)}}(a,b)$ is an injection
from $A\times B$ to $M\times M$ where $M$ is the set of midpoints of $A$
and $B$. This proves that $\#(A\times B)\leq \#(M\times M)$ as needed.

% As mentioned in the introduction, this construction extends to a higher
% number of sets involved: Given three sets $A,B,C$ in the cube, we 
% obtain injections from $A\times B\times C$ to $M\times M\times M$ where
% $M$ is the set of ``$1/3$-barycenters'' of $A$, $B$ and $C$, i.e.\ the
% points resulting from taking one third of the bits in each of $A$, $B$
% and $C$.

\bigskip

For later use, let us state a property of the coding
maps $\phi_c$ and $\Phi_c$. If $\Phi_c(a,b)=(m,m')$, we denote
$a=\phi_c^{-1}(m,m')$ and $b=\phi_{\bar c}^{-1}(m,m')=\phi_c^{-1}(m',m)$.

Let us equip the set of crossovers $C_r$ with
the distance
\[
d(c,c')\=\#(c\!\setminus\! c')+\#(c'\!\setminus\! c).
\]

\begin{Prop}[Decoding is isometric]
\label{prop:codingmetric}
%Let $a,b\in \{0,1\}^N$ and let $c,c' \in C_{d(a,b)}$. Then
%$d(c,c')=d(\phi_c(a,b),\phi_{c'}(a,b))$.
%
Let $m,m'\in \{0,1\}^N$. Let $c_1,c_2\in
C_{d(m,m')}$. Let $a_1=\phi^{-1}_{c_1}(m,m')$ and
$a_2=\phi^{-1}_{c_2}(m,m')$.
Then $d(a_1,a_2)=d(c_1,c_2)$.
\end{Prop}

\begin{proof}
Given $m$ and $m'$, modifying the crossover $c$ changes the preimage
$\phi^{-1}_c(m,m')$ by the same amount.
\end{proof}

\section{Concentration in the set of crossovers}

To get an improved inequality with positive curvature $K$, we will need
to study geometric properties of the set of crossovers; more precisely we
show that this set exhibits concentration of measure.  This is obtained
from the well-known concentration of measure in the permutation group by
a quotienting argument.  (We refer to \cite{Led01} for more
background about concentration of measure.) We first state concentration
in the permutation group under the form we need.

\begin{Lem}[Concentration in $S_n$]
Let $S_n$ be the permutation group on $\{1,\ldots,n\}$. Equip $S_n$ with
the distance $d(\sigma,\sigma')=\#\{i,\sigma(i)\neq \sigma'(i)\}$ for
$\sigma,\sigma'\in S_n$. Let $\nu$ be the uniform probability measure on
$S_n$.

Let $f:S_n\to\R$ be a $1$-Lipschitz function. Then $f$ satisfies the
concentration inequality
\[
\textstyle\nu(\{f\geq \int\! f\d\nu+t\})\leq e^{-t^2/2(n-1)}\qquad \forall
t\geq 0
\]
and the Laplace transform estimate
\[
\int\! e^{\lambda f}\d\nu\leq e^{\lambda \!\int\! f\d\nu \,+\,
(n-1)\lambda^2/2}\qquad
\forall \lambda\in \R.
\]
\end{Lem}

\begin{proof}
The second statement is Proposition~6.1 in~\cite{BHT06}.
The first statement follows by the
exponential Markov inequality.
\end{proof}

\begin{Prop}[The set of crossovers is concentrated]
\label{prop:crossconc}
Let $n\geq 1$ and let $C_n$ be the set of parts $c\subset\{1,\ldots,n\}$
with $\abs{\#c -n/2}<1$. Equip $C_n$ with the distance 
$d(c,c')\=\#(c\!\setminus\! c')+\#(c'\!\setminus\! c)$
as above
and with the uniform
probability measure $\mu$.

Let $f:C_n\to \R$ be a $1$-Lipschitz function. Then $f$ satisfies the
concentration inequality
\[
\textstyle\mu(\{f\geq \int\! f\d\mu+t\})\leq e^{-t^2/2n}\qquad \forall
t\geq 0
\]
and the Laplace transform estimate
\[
\int\! e^{\lambda f}\d\mu\leq e^{\lambda \!\int\! f\d\mu\, +\, n\lambda^2/2}\qquad
\forall \lambda\in \R.
\]
\end{Prop}

\begin{proof}
Let us begin with even $n$.
Then the natural
action of $S_n$ on $\{1,\ldots,n\}$ preserves $C_n$. Let us fix an origin
$c_0\=\{1,\ldots,n/2\}\in
C_n$ and define the projection map $\pi:S_n \to C_n$ by $\sigma\mapsto
\sigma(c_0)$. Each fiber of $\pi$ has the same cardinality
$\left((n/2)!\right)^2$. Moreover, if we equip $S_n$ and $C_n$ with the
distances as above,
then the map $\pi$ is
$1$-Lipschitz.

Thus, if $f:C_n\to \R$ is a $1$-Lipschitz function, the function
$\tilde f\= f\circ \pi$ is $1$-Lipschitz on $S_n$. So
$\tilde f$ satisfies the concentration
property $\nu(\{\tilde f\geq \int \tilde f\,d\nu+ t\})\leq
e^{-t^2/2(r-1)}$ where $\nu$ is
the uniform probability measure on $S_n$. Since all fibers of $\pi$ have
the same cardinality, $\pi$ sends $\nu$ to the uniform measure $\mu$ and
so the same estimate holds for $f$ in $C_n$ under $\mu$. The argument is identical for the Laplace transform
estimate.

For odd $n$ we proceed as follows. Let us fix
$c_0=\{1,\ldots,\lfloor n/2\rfloor\}\in
C_n$ and $c_1=\{1,\ldots,\lceil n/2\rceil\}\in
C_n$. Let us define the set $S^\ast_n \= S_n\!\times\!\{0\}\sqcup
S_n\!\times\!\{1\}$.
Define the map $\pi:S^\ast_n \to C_n$ by $(\sigma,i)\mapsto \sigma(c_i)$
for $i=0,1$. Then each fiber of $\pi$ has the same cardinality $\lfloor
n/2\rfloor!\,\lceil n/2\rceil!$. Let us equip $S^\ast_n$ with the
metric $d((\sigma,i),(\sigma',i'))=\abs{i-i'}+d(\sigma,\sigma')$.
Then one checks that $\pi$ is $1$-Lipschitz
from $S^\ast_n$ to $C_n$. (A more elegant construction would have used
$c\mapsto \bar{c}$ to get a group structure on $S_n^\ast$, but this has bad
metric properties.)

Given a $1$-Lipschitz function $f:C_n\to \R$, consider as above the
function $\tilde f\= f\circ \pi$ on $S_n^\ast$.
Applying, for instance, the technique of Theorem 4.2 in \cite{Led01} to get concentration of measure in
$S^\ast_n$ instead of $S_n$, we get that $\tilde f$ satisfies the Laplace
transform estimate
\[
\int e^{\lambda \tilde f}\,\d\nu\leq e^{\lambda \int \tilde f \d\nu +
(r-1)\lambda^2/2+\lambda^2/8}
\leq e^{\lambda \int \tilde f \d\nu +
r\lambda^2/2}
\]
 with $\nu$ the uniform
probability measure on $S_n^\ast$. This implies that $\nu(\{\tilde f\geq \int
\tilde f\,d\nu+ t\})\leq e^{-t^2/2r}$. Just as above, this estimate then holds for
$f$ on $C_n$.
\end{proof}

\begin{Cor}
\label{cor:crossoverseparation}
Let $A$ be a subset of the set of crossovers $C_n$ and let $\bar{A}\=\{\bar{c},\,c\in A\}$. Suppose that
$d(A,\bar{A})\geq k$. Then
\[
\#A\leq e^{-k^2/8n}\,\#C_n.
\]
\end{Cor}

\begin{proof}
Consider the function
$f:C_n\to \R$ given by $f(c)\=\frac12\left(d(c,A)-d(c,\bar{A})\right)$. This
function is $1$-Lipschitz, and takes values at least $k/2$ on $A$.
By symmetry the average of $f$ is $0$. So applying the
above, we get that the (relative) measure of $A$ in $C_n$ is at most
$e^{-k^2/8n}$.
\end{proof}

The following is a refined version of
Corollary¹\ref{cor:crossoverseparation}, in which the set $A$ is replaced
with a measure $\xi$, cardinals are replaced with entropies, and the
distance $d(A,\bar A)$ is replaced with $W_1(\xi,\bar \xi)$.

\begin{Cor}
\label{cor:crossoverw1sep}
Let $\xi$ be a probability measure on the set of crossovers $C_n$.
Let $\bar{\xi}$ be the complement of $\xi$ i.e.\ $\bar\xi(c)\=\xi(\bar
c)$ for $c\in C_n$. Then
\[
S(\xi)\leq \ln\#C_n-\frac{1}{8 n} W_1(\xi,\bar\xi)^2
\]
with $S$ the Shannon entropy.
\end{Cor}

\begin{proof}
The proof uses the following consequence of
Proposition~\ref{prop:crossconc}.

\begin{Lem}[$W_1H$ inequality for crossovers]
\label{lem:crossoverw1h}
Let $\xi$ be a probability measure on $C_n$. Then
\[
W_1(\xi,\mu)^2\leq 2n H(\xi|\mu)
\]
where $\mu$ is the uniform probability measure on $C_n$ and $H$ the
relative entropy.
\end{Lem}

Indeed, by a result of Bobkov and Götze (Theorem~3.1 in \cite{BG99}), the
inequality $W_1(\xi,\mu)^2\leq 2\gamma H(\xi|\nu)$ for all measures
$\xi$, is equivalent to the Laplace transform estimate
$\int\! e^{\lambda f}\d\mu\leq e^{\lambda \int\! f\d\mu +
\gamma\lambda^2/2}$ for all $\lambda\in \R$ and all $1$-Lipschitz functions
$f$. So the lemma is actually equivalent to
Proposition~\ref{prop:crossconc}.

Now, since $W_1(\xi,\bar\xi)\leq
W_1(\xi,\mu)+W_1(\mu,\bar\xi)=2W_1(\xi,\mu)$ by symmetry, we get
\[
H(\xi|\mu)\geq \frac{1}{8 n} W_1(\xi,\bar\xi)^2.
\]
Finally, using $H(\xi|\mu)=\ln\#C_n-S(\xi)$, this 
rewrites in terms of the Shannon entropy as 
\[
S(\xi)\leq \ln\#C_n-\frac{1}{8 n} W_1(\xi,\bar\xi)^2.
\]
\end{proof}

\section{Positively curved Brunn--Minkowski inequality}

Let us now prove Theorem~\ref{thm:BMcube}.
So let again $A,B$ be two nonempty subsets of the hypercube $X=\{0,1\}^N$,
and let $M$ be the set of midpoints of $A$ and $B$. We have to prove that
\[
\ln \#M\geq
\frac12\left(\ln\#{A}+\ln\#{B}\right)+\frac{K\,d(A,B)^2}{8},\qquad K =
\frac1{2N}.
\]

The difference with the case $K=0$ is that we now consider all crossovers
at once.
Let $C_r$ be the set of $r$-crossovers. Let
$Y\= \{(a,b,c), \, a\in A, \, b\in B, \, c\in C_{d(a,b)}\}$. 
Consider the map $f:(a,b,c)\mapsto \Phi_c(a,b)$ from $Y$ to
$M\times M$.
This map $f$ may not be one-to-one; but we will show that it is not
too-many-to-one. The idea is that, given a
pair of midpoints $(m,m')$, the geometry of $A$ and $B$ allows to guess,
to some extent, which crossover was used, so that the cardinality of
$f^{-1}(m,m')$ is bounded. (This is most clear when $A$ is a singleton
$\{00\ldots00\}$, in which case there is no ambiguity on the crossover:
every '1' in $m$ or $m'$ was taken from $B$.)

Let $Y_r\=\{(a,b,c)\in Y, \; d(a,b)=r\}$ and let likewise $(M\times
M)_r\=\{(m,m')\in M\times M, \; d(m,m')=r\}$. Now fix $(m,m')\in (M\times
M)_r$. The fiber $f^{-1}(m,m')$ is in bijection with the set $E$ of
crossovers $c\in C_r$ such that $\Phi_c^{-1}(m,m')\in A\times B$.
Consider, symmetrically, the set $E'=\{c\in C_r,\;
\Phi_c^{-1}(m,m')\in B\times A\}$. By definition 
$\Phi_c=(\phi_c,\phi_{\bar{c}})$, so the elements of $E'$ are
the complements of the elements of $E$.

We claim that $d(E,E')\geq d(A,B)$.  Indeed, 
 if $c\in E$,
$c'\in E'$ we have $\phi^{-1}_{c_1}(m,m')\in A$ and
$\phi^{-1}_{c'}(m,m')\in B$.
Since decoding is isometric (Proposition~\ref{prop:codingmetric})
we have $d(c,c') \geq d(A,B)$.

Corollary~\ref{cor:crossoverseparation} then states that
the cardinality of $E$ is at most $\#{C_r}
e^{-d(A,B)^2/8r}$. Since the cardinality of $E$ is also the
cardinality of the fiber $f^{-1}(m,m')$,
this
shows that the map $f:Y_r\to(M\times M)_r$ is at most $(\#{C_r}
e^{-d(A,B)^2/8r})$-to-one. Consequently, $\#{Y_r}\leq \#{C_r}
e^{-d(A,B)^2/8r} \,\#{(M\times M)_r}$.

Setting $(A\times B)_r\=\{(a,b)\in A\times B,\; d(a,b)=r\}$, 
we
have $\#{Y_r}=\#{(A\times B)_r}\times \#{C_r}$ so that
\[
\#{(M\times M)_r}\geq e^{d(A,B)^2/8r} \, \#{(A\times B)_r}.
\]
Finally, summing over $r$ from $1$ to $N$ we find
\[
\#(M\times M)\geq e^{d(A,B)^2/8N} \#(A\times B)
\]
which proves Theorem~\ref{thm:BMcube}.

\section{Entropy of the set of midpoints}

We now turn to the proof of Theorem~\ref{thm:entropicBM}.

Remember that, given $a$ and $b$ in the hypercube $X$, the midpoint
measure $\midm(a,b)$ is the uniform probability measure on all
midpoints of $a$ and $b$.
The
midpoint measure of two probability measures $\mu_A$ and $\mu_B$ is defined as
\[
\midm(\mu_A,\mu_B)\= \iint \midm(a,b)\,\d\mu_A(a)\d\mu_B(b)
\]
that is, the average of $\midm(a,b)$ where $a$ and $b$ are taken
\emph{independently} at random under $\mu_A$ and $\mu_B$.

The proof follows the same lines as in the deterministic case, using
probability measures instead of sets. The reader should think of the
probability measures below as being nothing but weighted sets, and their
Shannon entropy as being the logarithm of their cardinality. The main
differences are as follows:
\begin{itemize}
\item In the set-theoretic version, a key point was an estimation of
the cardinality of the fibers of the map $(a,b,c)\mapsto
(m,m')=\Phi_c(a,b)$. The lower bound on the cardinality of the set
$\{(m,m')\}$ followed. Here, we will use the associativity of Shannon
entropy to express the same relationship, yielding a lower bound on the
entropy of $(m,m')$ if the entropy of the fibers is known.
\item The final result involves $W_1(\mu_A,\mu_B)$ instead of $d(A,B)$.
In the set-theoretic version, we used the map $c\mapsto \bar c$ and the
fact that $\Phi_c(a,b)=\Phi_{\bar c}(b,a)$ to conclude that, if
$\Phi_c(a,b)=\Phi_{c'}(a',b')$ then $d(\bar c,c')=d(b,a')\geq d(A,B)$.
Then Corollary~\ref{cor:crossoverseparation} was used to bound the
cardinality of the set $E$ of such crossovers $c$ in a fiber.
The refined version uses the relation $d(\bar c,c')=d(b,a')$ to turn any
coupling between $E$ and $\bar E$, into a coupling between $A$ and $B$
with the same transportation distance. Then,
Corollary~\ref{cor:crossoverw1sep} is used as a refined version of
Corollary~\ref{cor:crossoverseparation} and yields a bound on the entropy
of the crossovers $c$ in a fiber.
\end{itemize}

So let $a$ and $b$ be independent random variables with law $\mu_A$ and
$\mu_B$. Let as above $C_r$ be the set of $r$-crossovers. Let $c$ be
a random variable uniformly distributed on $C_{d(a,b)}$, independent of $a$ and
$b$ conditionally to $d(a,b)$. Let us define the random variables
$m\=\phi_c(a,b)$ and $m'\=\phi_{\bar c}(a,b)$. Thus the law of $m$
is $\midm(\mu_A,\mu_B)$, as is the law of $m'$.

Let us slightly abuse notation and denote by $S((y))$ the Shannon entropy
of the law of a random variable $y$. We have $S((m,m'))\leq S((m))+S((m'))$
but since $m$ and $m'$ have the same law $\midm(\mu_A,\mu_B)$, we get
\[
S(\midm(\mu_A,\mu_B))\geq \frac12 S((m,m')).
\]

Consider as above the map $\Phi$ sending $(a,b,c)$ to $\Phi_c(a,b)=(m,m')$.
Let $Y_{(m,m')}$ be the law of $(a,b,c)$ knowing
$(m,m')$.
By the associativity of entropy, the Shannon entropy of the law of
$(m,m')$ is the
entropy of the law of $(a,b,c)$ minus the average entropy of
fibers of $\Phi$, namely:
\[
S((m,m'))=S((a,b,c))-\E S(Y_{(m,m')}).
\]

The first term is computed as follows. The random variables $a$ and $b$
are independent, and, conditionally to $d(a,b)$, the variable $c$ is
independent of $a$ and $b$ with law
the uniform distribution $U_{d(a,b)}$ on $C_{d(a,b)}$. So
\[
S((a,b,c))=S((a))+S((b))+\E S(U_{d(a,b)})=S(\mu_A)+S(\mu_B)+\E
\ln \# C_{d(a,b)}.
\]

Let us turn to the second term $\E S(Y_{(m,m')})$. This means
we have to evaluate the entropy of the fibers of $\Phi$, as in the
non-random case.

Let $E_{(m,m')}$ be the law of $c$ knowing $(m,m')$ (i.e., the
third marginal of $Y_{(m,m')}$).
Given $(m,m')$, the value of $c$ determines $a$ and $b$, and so,
$S((a,b,c)|(m,m'))=S((c)|(m,m'))$ i.e.\ 
\[
S(Y_{(m,m')})=S(E_{(m,m')})
\]
so that
\[
S((m,m'))=S(\mu_A)+S(\mu_B)+\E\ln\# C_{d(a,b)}-\E S(E_{(m,m')}).
\]

If, at this point, we apply the crude estimate $S(E_{(m,m')})\leq \ln
\#C_{d(m,m')}$, we get $S((m,m'))\geq S(\mu_A)+S(\mu_B)+\E
\ln \# C_{d(a,b)}-\E\ln\#  C_{d(m,m')}= S(\mu_A)+S(\mu_B)$ since
$d(a,b)=d(m,m')$. This implies $S((m))\geq \frac12 (S(\mu_A)+S(\mu_B))$
i.e.\ the case $K=0$ in the theorem.

As in the set-theoretic case, we will show that $E_{(m,m')}$ has small Shannon
entropy by using concentration properties in the set of crossovers. 
Corollary~\ref{cor:crossoverw1sep} tells us that
\[
S(E_{(m,m')})\leq \ln\#C_{d(m,m')}-\frac{1}{8 d(m,m')} W_1(E_{(m,m')},\bar E_{(m,m')})^2
\]
where $\bar E_{(m,m')}$ is the image of $E_{(m,m')}$ by $c\mapsto
\bar{c}$.
Thus,
we need to evaluate the distance between $E_{(m,m')}$ and
$\bar{E}_{(m,m')}$, as in the deterministic case.

Actually we only need an estimate on average over $(m,m')$. We claim that
\[
\E
W_1(E_{(m,m')},\bar{E}_{(m,m')})^2 \geq W_1(\mu_A,\mu_B)^2.
\]

Indeed, let us fix $(m,m')$ for now, and let $A_{(m,m')}$ and
$B_{(m,m')}$ be the laws of $a$ and $b$ knowing $(m,m')$, respectively.
Since $a=\phi_c^{-1}(m,m')$ and $b=\phi_{\bar{c}}^{-1}(m,m')$,
any
coupling between $E_{(m,m')}$ and $\bar{E}_{(m,m')}$ determines a coupling between
$A_{(m,m')}$ and $B_{(m,m')}$.
Moreover, since decoding is isometric by
Proposition~\ref{prop:codingmetric}, these
couplings will define the same transportation distance. So we get
$W_1(A_{(m,m')},B_{(m,m')})\leq W_1(E_{(m,m')},\bar{E}_{(m,m')})$. 

If for each $(m,m')$ we are given a coupling between $A_{(m,m')}$ and
$B_{(m,m')}$, by summation this defines a coupling between $\mu_A$
and $\mu_B$ and so
$W_1(\mu_A,\mu_B)\leq \E W_1(A_{(m,m')},B_{(m,m')})$. 
Thus 
$W_1(\mu_A,\mu_B)\leq
\E W_1(E_{(m,m')},\bar{E}_{(m,m')})$.
Then, by convexity
we get
\[
W_1(\mu_A,\mu_B)^2\leq\E W_1(E_{(m,m')},\bar{E}_{(m,m')})^2
\]
as announced.

Putting everything together and using that $d(m,m')=d(a,b)$,
we get
\begin{align*}
S((m,m'))
& = S((a,b,c))-\E S(Y_{(m,m')})
\\& = S(\mu_A)+S(\mu_B)+\E \ln\#C_{d(a,b)}-\E S(E_{(m,m')})
\\&\geq S(\mu_A)+S(\mu_B)+\E \ln\#C_{d(a,b)}-\E \ln\#C_{d(m,m')}+\E\left[
\frac{W_1(E_{(m,m')},\bar E_{(m,m')})^2}{8d(m,m')}\right]
\\&\geq
S(\mu_A)+S(\mu_B)+\frac{1}{8N}\,\E W_1(E_{(m,m')},\bar E_{(m,m')})^2
\\&\geq S(\mu_A)+S(\mu_B)+\frac{1}{8 N}W_1(\mu_A,\mu_B)^2
\end{align*}
and so
\[
S((m))\geq \frac12\left(S(\mu_A)+S(\mu_B)\right)+\frac{1}{16 N}W_1(\mu_A,\mu_B)^2
\]
which ends the proof.

%\signyo
%\signcv

\end{document}